\ifx\selectfont\undefined
\documentstyle[12pt]{article}
\else
\documentstyle[12pt,oldlfont]{article}
\fi

\makeatletter
\@addtoreset{equation}{section}
\makeatother

\newtheorem{fact}{Fact}[section]
\newtheorem{thm}[fact]{Theorem}
\newtheorem{prop}[fact]{Proposition}
\newtheorem{lemma}[fact]{Lemma}

\newcommand{\rem}{\noindent{\bf Remark{\ \ }}}
\newcommand{\exm}{\noindent{\bf Example{\ \ }}}
\newcommand{\proof}{{\noindent\bf Proof{\ \ }}}
\newcommand{\qed}{\bigskip\hfill\(\Box\)}
\newcommand{\Rn}[1]{\mbox{{\it I\kern -0.25emR}$^{\,{#1}}$}}
\newcommand{\Kn}[1]{\mbox{{\it I\kern -0.25emK}$^{\,{#1}}$}}
\newcommand{\RRn}[1]{\mbox{\scriptsize\Rn{#1}}}
\newcommand{\KKn}[1]{\mbox{\scriptsize\Kn{#1}}}
\newcommand{\Cn}[1]{\mbox{\it C\hspace*{-1.05ex}\rule{0.15ex}%
       {1.5ex}\hspace*{1.05ex}$^{\,{#1}}$}}
\newcommand{\NN}{\mbox{{\it I\kern -0.25emN}}}
\newcommand{\NNN}{\mbox{{\scriptsize\it I\kern -0.25emN}}}

\newcommand{\al}{\alpha}
\newcommand{\bt}{\beta}
\newcommand{\de}{\delta}
\newcommand{\ep}{\varepsilon}

\newcommand{\la}{\lambda}
\newcommand{\La}{\Lambda}
\newcommand{\spn}{\mathop{\rm span\,}}
\newcommand{\sgn}{\mathop{\rm sgn\,}}
\newcommand{\wt}{\widetilde}

\newcommand{\RRe}{\mathop{\rm Re\,}}
\newcommand{\IIm}{\mathop{\rm Im\,}}

\newcommand{\ie}{{\em i.e.,\/}}
\newcommand{\eg}{{\em e.g.,\/}}
\newcommand{\cf}{{\em cf.\/}}
\newcommand{\lsp}{\langle}
\newcommand{\rsp}{\rangle}

\begin{document}

\thispagestyle{empty}
\mbox{}
\vspace{1in}
\begin{center}
{\Large\bf Norms of Minimal Projections}
\end{center}
\vspace{15mm}
%
{\setbox2=\hbox{Mathematisches Seminar}
\parbox{\wd2}{%
  {Hermann K\"{o}nig}\footnotemark\hfill\\
  \box2\\
  Universit\"{a}t Kiel\hfill\\
  Kiel, {Germany}\hfill}
\footnotetext{${}^{,2}\,$During the work on this paper both authors
  were partially  supported by NATO Collaborative
Research Grant CRG 920047.}
\ {\hfill} \
\setbox2=\hbox{{Nicole Tomczak-Jaegermann}\footnotemark}
\parbox{\wd2}{%
  \box2
  Department of Mathematics\hfill\\
  University of Alberta\hfill\\
  Edmonton, Alberta, {Canada}\hfill}
}
\vfill

\begin{abstract}
\noindent
It is proved that the projection constants of two- and
three-di\-men\-sio\-nal spaces are bounded by $4/3$ and $(1+\sqrt 5)/2$,
respectively. These bounds are attained precisely by the spaces
whose unit balls are the regular hexagon and dodecahedron.
In fact,  a general inequality for the projection constant
of a real  or  complex $n$-dimensional space is obtained
and the question of equality therein is discussed.
\end{abstract}

\vspace{2cm}
\clearpage


\section{Introduction and the main results}

In this paper we prove results on upper estimates for the norms
of minimal projections onto finite-dimensional subspaces of Banach spaces,
which are optimal in general.
By Kadec--Snobar [KS], onto $n$-dimensional spaces there are
always projections of norm smaller than or equal to $\sqrt n$.
General bounds for these so-called projection constants
were further studied
by various authors, including
Chalmers, Garling, Gordon, Gr\"{u}nbaum,
K\"{o}nig, Lewis and Tomczak-Jaegermann ([GG], [G], [KLL], [KT],
[L], [T]).
Some other aspects of minimal projections, like the existence
or norm estimates for concrete spaces, were investigated
by many authors, among them \eg\ Chalmers, Cheney, Franchetti
([CP], [IS], [FV]).

In [KT] a very tight formula for the projection constants of spaces
with enough symmetries was shown. We now prove that this formula holds
for arbitrary spaces, and study cases of equality.

The formula yields, in particular, that the projection constant
of any real (resp.\ complex) 2-dimensional space is bounded by
$4/3$ (resp. $(1+\sqrt 3)/2$). Up to isometry, there is just one
space (in each case) attaining the bound. The values for
3-dimensional spaces are $(1+ \sqrt 5)/2$ (resp. $5/3$).
In the real case, the unique extremal spaces are those whose unit
balls are the regular hexagon and the regular dodecahedron.
The $4/3$-result solves a problem   of Gr\"{u}nbaum [G].
A proof of this fact has also been announced by Chalmers et al. [CMSS];
it is our understanding that their argument is incomplete
as of now.

The authors would like to thank to J.~J.~Seidel
for valuable remarks concerning equiangular lines.

\bigskip

We use standard Banach space notation, see \eg\  [T.2].
By $\Kn{}$ we denote the scalar field, either the real numbers
$\Rn{}$ or the complex numbers $\Cn{}$.
The {\em relative projection constant} of a (closed) subspace
$E$ of a Banach space $X$ is defined by
$$
\la (E, X):= \{ \|P\| \, \mid \, P: X \to E \subset X
\mbox{\ is a linear projection onto\ }E \},
$$
the {\em (absolute) projection constant} of $E$ is given by
\begin{equation}
\la(E):= \{ \la(E, X)  \mid  X
\mbox{\, is a Banach space containing\, }
E \mbox{\, as a subspace}\}.
  \label{ab_proj}
\end{equation}

Any separable Banach space $E$ can be embedded isometrically
into $l_\infty$.  For any such embedding,
$\la (E) = \la (E, l_\infty)$, \ie\ the supremum in   (\ref{ab_proj})
is attained. We can therefore restrict our attention
to  finite-dimensional subspaces
$E \subset l_\infty$.
Also note that $\la (l_\infty^n)= \la (l_\infty)=1$.

Let $n \in \NN$ be a positive integer,
$\lsp \cdot, \cdot \rsp$ denote the standard scalar  product
in $\Kn{n}$ and let $\|\cdot\|_2 = \sqrt {\lsp \cdot, \cdot \rsp}$.
For $ N \in \NN$,
vectors $x_1, \ldots, x_N \in \Kn{n}$ spanning lines
in $\Kn{n}$ are called {\em equiangular} provided that
there is $0 \le \al <1$ such that
$$
\|x_i\|_2 =1 \mbox{\ and \ } |\lsp x_i, x_j \rsp| = \al
\mbox{\ for\ }i \ne j,\  i,j= 1, \ldots, N.
$$

Put
\begin{equation}
N(n):=\left\{
\begin{array}{ll}
 n (n+1)/2 & \mbox{if  $\Kn{} = \Rn{}$}\\
  n^ 2 & \mbox{if  $\Kn{} = \Cn{}$.}
\end{array}
\right.
\label{Nn}
\end{equation}
By Lemmens--Seidel [LS] and Gerzon, in $\Kn{n}$
there are at most $N(n)$ equiangular vectors.
(Indeed, the hermitian rank 1 operators $x_i \otimes x_i$
are linearly independent
in a suitable  real linear space of operators.)
This bound is attained for $n = 2, 3, 7, 23$
if $\Kn{} = \Rn{}$ and  for  $n = 2, 3$ if
$\Kn{} = \Cn{}$. If the bound is attained, necessarily
$\al = 1/ \sqrt {n+2}$ if $\Kn{} = \Rn{}$
and  $\al = 1/ \sqrt {n+1}$ if $\Kn{} = \Cn{}$.

Our main result is
\begin{thm}
  \label{main}
{\bf (a)}
The projection constant of any $n$-dimensional
normed space $E_n$ is bounded by
{\renewcommand{\arraystretch}{1.5}
\begin{equation}
  \label{ineq_main}
\la (E_n) \le  \left\{
         \begin{array}{ll}
\bigl(2+ (n-1)\sqrt{n+2}\bigl)/ (n+1) &
     \mbox{in the real case,}\\
\bigl(1+ (n-1)\sqrt{n+1}\bigl)/ n  &
     \mbox{in the complex case.}
        \end{array}
\right.
\end{equation}
}

{\bf (b)}
Given $\Kn{}$ and $n \in \NN$, there  exist $n$-dimensional spaces
$E_n$ for which the bound is attained if and only if
there exist $N(n)$ equiangular vectors in $\Kn{n}$. In this case, such
a space $E_n$ can be realized as an isometric subspace of
$l_\infty^{N(n)}$, and the orthogonal projection is a minimal
projection onto $E_n$.

{\bf (c)}
For $\Kn{}= \Rn{}$ and $n=2, 3, 7, 23$, there are {\rm unique}
spaces $E_n$ (up to isometry)  attaining the bound (\ref{ineq_main});
for $\Kn{}= \Cn{}$ and $n=2, 3$ such spaces also exist.
For $\Kn{}= \Rn{}$ and $n=2, 3$ the unit balls of $E_n$
are the regular hexagon and the regular dodecahedron, respectively.
\end{thm}

\noindent{\bf Remarks\ }
{\bf (i)}
The right hand side of (\ref{ineq_main}) equals the
bound  $f(n, N(n))$ derived in [KLL]
for the relative projection constant
of an $n$-dimensional space in an $N(n)$-dimensional superspace.

{\bf (ii)}
The bounds in (\ref{ineq_main}) are  of the order
$\sqrt{ n} - 1/ \sqrt{ n} + 2/n$ if $\Kn{} = \Rn{}$
and  $\sqrt{ n} - 1/2 \sqrt{ n} + 1/n$ if $\Kn{} = \Cn{}$,
for large $n \in \NN$.

\bigskip

To prove just (\ref{ineq_main}) it would
suffice (by approximation) to consider polyhedral spaces
$E \subset l_\infty^N$ for an arbitrary $N \in \NN$
(the ``finite'' case); in which case the
proofs of most of the results which follow can be simplified.
For the examination of the equality in (\ref{ineq_main}) and
the uniqueness we need,
however, the general (``infinite'')
case of  $E \subset l_\infty$ as well,
even though the spaces attaining the bound
(\ref{ineq_main}) turn out in the end to be
polyhedral.

To unify the notation in the finite and infinite case which we
would like to discuss simultaneously, we set $T= \{1, \ldots, N\}$,
for some $N \in \NN$,
in the finite case and  $T = \NN$, in the infinite case.
In particular, $l_\infty(T)$ denotes $l_\infty^N$ in the
former case and $l_\infty$ in the latter case.

If $\mu= (\mu_t)_{t \in T}$ is a probablity measure on $T$,
and $ 1 \le p <\infty$, we let
$$
l_p(T, \mu):= \{(\xi_t)_{t\in T}\, \mid \,
\|(\xi_t)\|_{p,\mu}= \bigl(\sum_{t \in T}
         |\xi_t|^p \mu_t \bigr)^{1/p}<\infty\}.
$$
For  a subspace  $E \subset l_\infty (T)$,
we denote by  $E_{p, \mu}$
the same space $E$ considered as a subspace of $l_p(T,\mu)$,
via the embedding $l_\infty(T) \to l_p(T,\mu)$.

Finally, for $N \in \NN$,
by  $R_N: l_p \to l_p^N$ we denote the projection onto
the first $N$ coordinates,
acting in an appropriate sequence space
($1 \le p \le \infty$).

\bigskip

Let $n \in \NN$. The set ${\cal F}_n$ of all $n$-dimensional spaces,
equipped with the (logarithm of the) Banach--Mazur distance,
is a compact metric space, \cf\ \eg\ [T.2]. The projection
constant $\la$, as a  function
$\la: {\cal F}_n \to \Rn{}^+$, is continuous with
respect to this metric, and hence the supremum
$\sup_{E \in {\cal F}_n} \la(E)$ is attained:
there is $F \in {\cal F}_n$ with
\begin{equation}
\label{max_proj}
\la (F) = \sup \{ \la(E)\,|\, {E \in {\cal F}_n}\}.
\end{equation}

The proof of the bound (\ref{ineq_main}) is based upon
an estimate in terms of orthonormal systems,
which, in fact, is a characterization of the maximal
projection constant, and it
seems to be of independent interest.

\begin{thm}
\label{la=sup}
Let $n \in \NN$. Then
\begin{equation}
\label{la=sup_eq}
\max_{E \in {\cal F}_n} \la(E)
=  \sup_{\mu} \sup_{\{f_j\}}
\sum_{s, t\in \NNN} | \sum_{j=1}^n f_j(s)
           \overline{f_j (t)}| \mu_s \mu_t,
\end{equation}
where the outside supremum runs over  the
set  of all discrete probability measures $\mu= (\mu_t)_{t}$
on $\NN$ and the
inside supremum runs over all
orthonormal systems $\{f_j\}$ in $ l_2(\NN,\mu)$.
The double supremum in (\ref{la=sup_eq}) is attained
for some $\mu$ and
$\{f_j\}\subset l_2(\NN,\mu)\cap l_\infty$.
In this case,
the space $E = \spn \{f_1, \ldots, f_n\} \subset l_\infty$
has  maximal projection constant.
The square function
$(\sum_{j=1}^n |f_j(s)|^2)^{1/2}$
is constant $\mu$-a.e.\ in the extremal case.
\end{thm}

In the extremal case the support of $\mu$ can  be finite;
and in dimensions $n = 2, 3$ it is  actually so. The upper
estimate in  (\ref{la=sup_eq}) relies on an idea of
Lewis [L]. To prove Theorem~\ref{main}, we then have to find
an upper estimate for the right hand side of
(\ref{la=sup_eq}). In certain dimensions ($n=2, 3$,
and in the real case additionally  $n=7,23$),
we find the  exact value of  (\ref{la=sup_eq});
for other $n \in \NN$, the expression in
(\ref{la=sup_eq}) might be possibly used to slightly
improve (\ref{ineq_main}).

\section{Projection constants and trace duality}

In this section, we prove the upper bound for
$\max \{\la(E)\, \mid\, E \in {\cal F}_n\}$
in (\ref{la=sup_eq}). The argument is based
on trace duality.
For the convenience of the general reader, we try to use only
basic Banach space theory.
The first lemma is similar to Lemma 1
of [KLL].

\begin{lemma}
  \label{dual}
Let $E \subset l_\infty(T)$ be  a finite-dimensional subspace,
where $T = \{1, \ldots, N\}$, or  $T = \NN$.
There exists a map $u: l_\infty(T) \to l_\infty(T)$ with $u(E) \subset E$
such that
$$
\la (E) = \mathop{\mbox{tr\,}} (u: E \to E) \qquad \mbox{and} \qquad
\sum_{t \in T} \|u e_t\|_\infty =1.
$$
Here $(e_t)_{t \in T}$ denotes the standard unit vector basis
in $ l_\infty(T)$.
\end{lemma}

In fact, for any map $u$ with $u(E) \subset E$ and
$\sum_{t \in T} \|u e_t\|_\infty =1$ one has
$\mbox{tr\,} (u: E \to E) \le \la (E)$,
see (\ref{six}) below.

\medskip

\noindent\proof
Since $E$ is finite-dimensional, there exists a minimal  projection
onto $E$, say $P_0: l_\infty (T) \to  E \subset l_\infty (T) $ with
$\|P_0\| = \la = \la (E) <\infty$ (\cf\ [BC], [IS]).
Let ${\cal F}(l_\infty, l_\infty)$ denote the space
of finite-rank operators on $l_\infty =l_\infty(T)$,
equipped with the operator norm. The sets
$$
A= \{ S \in {\cal F}(l_\infty, l_\infty) \, \mid \, \|S\| < \la \}
$$
and
\begin{eqnarray*}
B &=& \{ P \in {\cal F}(l_\infty, l_\infty) \, \mid \,
P = P_0 + \sum_{i=1}^m x_i^* \otimes x_i \\
&&\qquad\qquad \mbox{for some \ } x_1, \ldots, x_m \in E,\
      x_1^*, \ldots, x_m^* \in E^{\bot}\subset l_\infty^*,\ m \in \NN \}
\end{eqnarray*}
are convex and  disjoint, since $B$ consists of projections onto $E$
and $\|P\|\ge \la$ for every projection $P$.
Since $A$ is open, by the Hahn--Banach theorem there is a functional
$\varphi \in {\cal F}(l_\infty, l_\infty)^*$ of norm
$\|\varphi \| =1$ such that $\varphi (P_0) \in \Rn{}$ and
for $S \in A$ and $P \in B$ we have
$$
\RRe \varphi (S) <\la \le \RRe \varphi (P).
$$

By  the trace duality, $\varphi$ is represented by a map $v$
defined on $l_\infty(T)$.
In the case  $T = \{1, \ldots, N\}$, the operator norm
of $w \in {\cal F}(l_\infty, l_\infty)$
is just
$\sup_{t \in T} \|w^* e_t\|_1$,
so the dual norm is
$\sum_{t \in T} \|v e_t\|_\infty$.
If $T = \NN$, define a linear operator
$v: l_\infty \to l_\infty^{**}$
by
$\lsp v(x), x^*\rsp = \varphi ( x^* \otimes x)$
for $x \in l_\infty$ and $x^* \in l_\infty^*$.
Writing any $S \in {\cal F}(l_\infty, l_\infty)$
as
$S= \sum_{i=1}^m x^* \otimes x_i$, one finds that
$\varphi (S) = \mbox{tr\,}(v\,S)$, with the {\em integral norm}
$i(v)$ equal to
$$
i(v) = \sup_{S \in {\cal F}(l_\infty, l_\infty)}
\mbox{tr\,}(v\,S) / \|S\| =1.
$$

Let $x^* \in E^{\bot}$, $x \in E$. Then
$\la \le \RRe \varphi (P_0 + x^* \otimes x) = \la +
\RRe \mbox{tr\,}(v\,(x^* \otimes x))$.
Hence $\RRe \lsp vx, x^*\rsp \ge 0$ for all
$x^* \in E^{\bot}$, $x \in E$, which implies
$\lsp vx, x^*\rsp =0$. Thus
$v(E)\subset E^{\bot \bot}=E \subset l_\infty$,
in view of $\dim E < \infty$.
Let $Q: l_\infty^{**} \to l_\infty$ be the canonical projection
onto $l_\infty$ with $\|Q\| =1$.
Let $u:= Q\,v: l_\infty \to l_\infty$.
Then $u(E) \subset E$ and, since $Qx=x$ for $x \in E$,
we have $u\, P_0 = v\, P_0$. Furthermore,
$i(u) \le \|Q\| i(v) =1$ and
$$
\la(E) = \la = \varphi (P_0)=
\mbox{tr\,}(u\,P_0) = \mbox{tr\,}(u: E \to E).
$$

Let $R_N:  l_\infty \to l_\infty^N$  be the natural projection
and  let $u_N:= R_N u:l_\infty \to l_\infty^N$.
Then $i(u_N ) \le 1$ and, similar as in the case of
$T = \{1, \ldots, N\}$ discussed above, this norm,
being dual to the operator norm on
${\cal F}(l_\infty^N , l_\infty)$, is equal to
$i(u_N ) = \sum_{i \in T} \|u_N  e_t\|_\infty \le 1$.
Taking the limit as $N  \to \infty$
(first for finite sums in $t$) we get that
$ \sum_{i \in T} \|u e_t\|_\infty \le 1$.
In fact we have the equality.
\qed

The following upper estimate is a consequence of Lemma~\ref{dual}
and relies essentially on an idea of Lewis [L].

\begin{prop}
  \label{lewis}
Let $E \subset l_\infty(T)$ be an $n$-dimensional subspace,
where $T = \{1, \ldots, N\}$, or  $T = \NN$.
There is a discrete probability measure $\mu= (\mu_t)_{t \in T}$
on $T$, $\|\mu\|_1 =1$ such that for any orthonormal basis
$(f_j)_{j=1}^n$ in $E_{2,\mu}$ we have
$$
\la (E) \le \sum_{s, t\in T} | \sum_{j=1}^n f_j(s)
                   \overline{f_j (t)}| \mu_s \mu_t.
$$
\end{prop}

Note that the double sum in this proposition is  finite,
since $f_j \in l_2(T, \mu)$.

\medskip

\noindent \proof
Let $u: l_\infty(T) \to  l_\infty(T) $ be as in Lemma~\ref{dual}
and put $\mu_t = \|ue_t\|_\infty$, for $t \in T$.
Then $\mu$ is a probability measure on $T$, $\sum_{t \in T} \mu_t =1$.
For every $N \in \NN$, let $u_N= R_N u: l_\infty \to l_\infty^N$.
Then we have
$\|u_N: l_1(T,\mu) \to l_\infty(T)\| \le 1$.
Indeed,  the extreme points of the unit ball in $ l_1(T,\mu)$
are, up to a multiple of modulus 1,
of the form $e_t/\mu_t$ for $t \in T$ and we have
$\|u_N(e_t/\mu_t)\|_\infty = \|u_N(e_t)\|_\infty/\mu_t \le 1$.

Now consider $E_{2, \mu} \subset l_2(T, \mu)$ and
fix an arbitrary orthonormal  basis $\{f_j\}_{j=1}^n$
in $E_{2, \mu}$. Then
\begin{eqnarray*}
\la (E) &=& \mbox{tr\,}(u: E \to E)
    = \mbox{tr\,}(u: E_{2,\mu} \to E_{2,\mu})\\
&=& \sum_ {j=1}^n \lsp u f_j, f_j\rsp_{l_2(\mu)}
= \lim _{N\to \infty} \sum_ {j=1}^n \lsp u_N f_j, f_j\rsp_{l_2(\mu)}.
\end{eqnarray*}
The second equality is purely algebraic; for the last one  use the fact
that $\lsp R_N g, h\rsp $ tends to $\lsp  g, h\rsp $ as $N \to \infty$, for
all $g, h \in l_2(T, \mu)$. Hence, by Lewis' idea
of  how to use the
bound for the norm of $u_N$ considered above, we have
\begin{eqnarray*}
  \la(E) &\le& \limsup_{N \to \infty}
            \sum_{ t\in T} | \sum_{j=1}^n u_N f_j(t)
                   \overline{f_j (t)}| \mu_t \\
 &\le& \limsup_{N \to \infty}
            \sum_{ t\in T} \|u_N \bigl(\sum_{j=1}^n
                   \overline{f_j (t)} f_j\bigr)\|_\infty \mu_t \\
 &\le& \sum_{s, t\in T} | \sum_{j=1}^n f_j(s)
                   \overline{f_j (t)}| \mu_s \mu_t,
\end{eqnarray*}
as required.
\qed

\section{Square function in an extremal case}

As a consequence of Proposition~\ref{lewis},
given a space $E \subset l_\infty(T)$,
an upper bound  for $\la(E)$ would
follow from an upper estimate for  the quantity
\begin{equation}
\phi(n, T) = \sup_{\mu \in {\cal M}} \sup_{\{f_j\}}
\sum_{s, t\in T} | \sum_{j=1}^n f_j(s) \overline{f_j (t)}| \mu_s \mu_t,
  \label{sups}
\end{equation}
where the outside supremum runs over  the set
${\cal M}$
of all discrete
probability measures $\mu$ on $T$ and the
inside supremum runs over all
orthonormal bases $\{f_j\}$ in $E_{2,\mu} \subset l_2(T,\mu)$.

To estimate (\ref{sups}), we first show,  using Lagrange multipliers,
that the {\em  square function}
of an extremal system $\{f_j^0\}$
is constant $\mu$-a.e.

We will be mainly concerned with the situation when $\phi$ really
increases at the dimension $n$,
\begin{equation}
  \label{phi}
  \phi(n_1, T_1) < \phi(n, T) \quad {\rm whenever\ } n_1 < n \ {\rm
    and\ } T_1 \subset T.
\end{equation}

\begin{prop}
  \label{sq_funct}
  Let $n \in \NN$  and  let  $T = \{1, \ldots, N\}$, or  $T = \NN$
  satisfy
  (\ref{phi}).  Assume that $\mu^0 \in {\cal M}$ and an
  orthonormal system
  $\{f_j^0\}_{j=1}^n $ in $l_2(T,\mu)$ attains the supremum
\begin{equation}
    \sum_{s, t \in T} | \sum_{j=1}^n f_j^0(s) \overline{f_j^0 (t)}|
    \mu_s^0  \mu_t^0
 =   \phi(n, T).
    \label{extr}
  \end{equation}
  Then the square function $f^0$ is constant $\mu$-a.e.,
$$
f^0(s):= \Bigl(\sum_{j=1}^n |f_j^0(s)|^2\Bigr)^{1/2} =
 \left\{
\begin{array}{ll}
  \sqrt{n} & \mbox{if $\mu_s^0 \ne 0$}\\ \ 0 & \mbox{if $\mu_s^0 = 0$}
\end{array}
\right.
$$
\end{prop}

\medskip
First notice that if $\mu_s^0 = 0$ for some $ s \in T$ then
$f_j^0(s)=0$ for $j=1, \ldots, n$, hence also $f^0(s)=0$.
Indeed, otherwise decreasing $|f^0_j(s)|$ would allow us to multiply
all the remaining $|f^0_j(t)|$ for $t \ne s$, by $\xi >1$, thus increasing the v
of the sum in (\ref{extr}).

Condition (\ref{phi}) implies that the matrix $(f_j^0(s))$
does not split into a non-trivial block diagonal sum of smaller
submatrices.

\begin{lemma}
  \label{chain}
  For all $ l, m= 1,\ldots,n$ we have
  \begin{equation}
\exists  l = l_0, \ldots, l_\rho=m
\ \forall 1 \le r \le \rho
\ \exists s \in T, \mu^0_s \ne 0
\qquad f^0_{l_{r-1}}(s)  f^0_{l_r}(s)  \ne 0.
    \label{product}
  \end{equation}
\end{lemma}

\proof
For $0 < \tau \le 1$, by ${\cal M}_\tau$ denote
the set of all discrete measures $\mu$ on $T$ such that
$\mu (T) = \tau$. By  $ \phi(n, T,\tau)$ denote the corresponding
supremum, analogous to (\ref{sups}), so that
$\phi(n, T)=  \phi(n, T,1)$.

It is easy to check that
$\phi(n, T, \tau)= \tau \phi(n, T, 1)$.
Moreover,
$\phi(n_1, T_1, 1) \le \phi(n, T, 1)$
if $n_1 \le n$ and $T_1 \subset T$.

Let $J_1\subset \{1, \ldots, n\}$ be a maximal set such
that (\ref{product}) is
satisfied for all $l, m \in J_1$
and let $J_2 = \{1, \ldots, n\} \backslash J_1$ be
the complement of $J_1$. Clearly, $J_1 $ is non-empty.
Let $T_1\subset T$ be the set of all $s$ such that
$f^0_{ j}(s) \ne 0$ for some $j \in J_1$, let
$T_2 =T \backslash T_1$.  By the maximality of $J_1$
and the definition of $T_1$ we have
$$
f^0_{ j}(s) =0 \quad {\rm whenever}\quad
(s,  j) \in (T_2 \times J_1) \cup (T_1 \times J_2).
$$

Denote by $\Phi$ the function whose supremum is taken in (\ref{sups}),
and by  $\Phi_1 $ and $\Phi_2$  the functions given by
the analogous formulas,
with the summation extended over
$s, t \in T_1$ and $j \in J_1$ for $\Phi_1$,   and over
$s, t \in T_2$ and $j \in J_2$ for $\Phi_2$.
We have $\Phi = \Phi_1 + \Phi_2$.
Moreover, as the functions $\Phi_i$ involve only sets $J_i$
and $T_i$, then
$\Phi_i(z_{s j} , \la_s ) \le \phi(n_i, T_i, \tau_i)$,
where $n_i = |J_i|$ and $\tau_i = \sum_{s \in T_i}
\la_s^2$, for $i=1, 2$. Thus
\begin{eqnarray}
\label{convex}
\lefteqn{\phi(n, T, 1) = F(z_{s j} , \la_s )= F_1(z_{s j} , \la_s ) +
  F_2(z_{s j} , \la_s )}\\ &&\le  \phi(n_1, T_1, \tau_1)+ \phi(n_2,
T_2, \tau_2) = \tau_1 \phi(n_1, T_1,1) + \tau_2 \phi(n_2,
T_2,1).\nonumber
\end{eqnarray}
Since $\tau_1 + \tau_2 =1$ and $n_1 >0$, the assumption  (\ref{phi})
implies that the  inequality in (\ref{convex}) is not
possible unless $\tau_2=0$ and  $n_1 =n$.
Thus $J_1 = \{1, \ldots, n\}$ and hence
(\ref{product}) holds for all $l$ and  $m$, as required.
\qed

 To simplify the orthogonality conditions, we let
\begin{equation}
  Z_{sj} = f_j(s) \sqrt{\mu_s} \quad {\rm and}\quad \La_s =
  \sqrt{\mu_s} \quad {\rm for\ } s \in T,\, j = 1, \ldots, n.
  \label{subst}
\end{equation}

Given the matrix $(Z_{s j})_{s\in T,1 \le j\le n}$
we consider
``short'' vectors $Z_s = (Z_{s j})_{j} \in \Kn{n}$,
$s \in T$,
and ``long'' vectors $\wt{Z}_j = (Z_{s j})_{s\in T} \in l_2$,
$j = 1, \ldots, n$.
The natural scalar product both in $\Kn{n}$ and in $l_2$
will be denoted by $\lsp \cdot, \cdot \rsp$.

We work with  the function $  F(Z_{s j}, \La_s) $
defined by
\begin{equation}
  F(Z_{s j}, \La_s) = \sum _{s,t \in T} |\lsp Z_s, Z_t \rsp| \La_s \La_t
  = \sum _{s,t\in T} | \sum _{j=1}^n Z_{s j} \overline{Z}_{t j}| \La_s
  \La_t.
        \label{obj_funct}
\end{equation}

\medskip
\noindent{\bf Proof of Proposition~\ref{sq_funct}}
{\bf (a)}
First let $\Kn{}=\Rn{}$ and $T = \{1, \ldots, N\}$.
We use Lagrange multipliers.
Clearly, the supremum $\phi (n,N)$ described in
(\ref{extr}) is equal to the
maximum of $F$ on the surface given by the  conditions
\begin{eqnarray}
  \label{constr_1}
  G_{l m}(Z_{s j}, \La_s) &:=& \lsp\wt{Z}_l, \wt{Z}_m\rsp - \de_{l m} =
  \sum _{s\in T} Z_{s l} \overline{Z}_{s m} - \de_{l m} = 0 \nonumber\\
  &&\quad\quad \quad {\rm for}\ 1 \le l \le m \le n \\
  G_0 (Z_{s j}, \La_s) &:=& \lsp \wt \La, \wt \La \rsp -1
       = \sum_{s\in T} \La _s ^2 -   1 =0
                \quad {\rm for}\ s \in T.
  \label{constr_2}
\end{eqnarray}
%
The supremum is attained for  a sequence of non-negative
$\La_s$;  if we set  $z_{sj} := f^0_j(s) \sqrt{\mu^0_s}$
and $\la_s := \sqrt{\mu^0_s}$ for $s \in T$,
$ j = 1, \ldots, n$, then
$F$ attains  its maximum at $ (z_{s j}, \la_s)$.

Consider the  Lagrange function $L$ defined by
$$
2 L(Z_{s j}, \La_s) = F(Z_{s j}, \La_s) -
\sum_{l \le m} \wt{\gamma}_{l m} G_{l m} (Z_{s j}, \La_s) -
\bt G_0 (Z_{s j}, \La_s).
$$

Assume that  $ (z_{s j}, \la_s)$ is  a point
where $F$ attains a  local maximum
subject to (\ref{constr_1}) and (\ref{constr_2}).
If for some $1 \le s, t  \le N$ we had
$\lsp z_s , z_t  \rsp =0$, we would leave this term
out from the sum defining $F$.
This would lead to a new sum,
defining the new function $F_1$.
Clearly, $F_1 \le F$ and $\max F_1 = \max F$.
The maximum is  attained at the same point $(z_{s j} , \la_s )$ and
the function $F_1$  is $C^2$ in the neighborhood of this point.
Moreover, by setting $\sgn 0 = 0$,
in the formulas for derivatives which follow
we will still be able to extend the sums over all indices $s, t$.

To use standard necessary conditions for Lagrange multipliers
we first check that the point   $ (z_{s j}, \la_s)$ is
regular. This means that the gradients
$\nabla G_{0}$ and $\nabla G_{lm}$
for $1 \le l \le m \le n $, are linearly
independent vectors in  $\Rn{N(n+1)}$.

Denoting vectors $(z_{sj})_{s=1}^N$ by  $\wt{z}_{j}$
for $j=1, \ldots, N$ and
$(\la_{s})_{s=1}^N$ by  $\wt{\la}$,
by a straigtforward differentiation
with respect to $Z_{sj}$ and $\La_s$ we get,
for $1 \le l, m \le n$ and $l < m$,
\begin{equation}
\nabla G_{0}
=  \left( \begin{array}{c}
0\\\vdots\\ 0\\ 2\wt{\la}
\end{array}\right),
\quad
\nabla G_{lm}
=  \left( \begin{array}{c}
0\\\vdots\\ \wt{z}_{m}\\\vdots\\
\wt{z}_l\\ \vdots\\ 0\\ 0
\end{array}\right),
\quad
\nabla G_{mm}
=  \left( \begin{array}{c}
0\\\vdots\\ 2 \wt{z}_{m}
\\ \vdots\\ 0\\ 0
\end{array}\right).
  \label{gradients}
\end{equation}
In the formula for $\nabla G_{lm}$, with $l <m$,
$\wt{z}_{m}$ stays on the $l$th place and
$\wt{z}_{l}$ stays on the $m$th place; and
in the formula for $\nabla G_{mm}$,
$2 \wt{z}_{m}$ stays on the $m$ place.
Since $\la_s \ne 0$ for $s=1, \ldots, N$,
the linear independence of the gradient vectors
(\ref{gradients}) follows directly from the
linear independence of the vectors $\wt{z}_{l} \in  \Rn{N}$,
for $l=1, \ldots, n$; the latter fact is an immediate consequence
of the orthogonality, hence linear independence,
of the system $\{f^0_l\}_{l=1}^n$.

Now, the first order condition
for Lagrange multipliers
states that
there exist multipliers
$\wt{\gamma}_{l m}$ and $\bt$
such that after setting
$$
\gamma_{l m} =   \frac{1}{2}  \left\{
\begin{array}{ll}
                        \wt{\gamma}_{l m}   & \mbox{if  $l< m$}\\
                        \wt{\gamma}_{ m l } & \mbox{if  $m< l$}\\
                      2 \wt{\gamma}_{ l l } & \mbox{if  $m= l$}
\end{array}
 \right.
$$
we have
\begin{eqnarray}
\label{dl/dz}
\frac{\partial L}{\partial  Z_{s l}} &=&
\sum _{t\in T} \sgn \lsp z_s, z_t \rsp  z_{t l} \la_s \la_t
- \sum_{m=1}^n {\gamma}_{l m} z_{s m} =0 \nonumber\\
&&\qquad\qquad  {\rm for}\ s \in T, l=1, \ldots, n \\
\label{dl/dla}
\frac{\partial L}{\partial  \La_{s}} &=&
\sum _{t\in T} |\lsp z_s, z_t \rsp| \la_t
- \bt  \la_s =0 \quad  {\rm for}\ s \in T.
\end{eqnarray}

First we simplify (\ref{dl/dz})
by a suitable orthogonal transformation.
Define two $N \times N$ matrices $A$ and $B$ by
\begin{equation}
\label{a_b}
A = \left(\sgn \lsp z_s, z_t \rsp  \la_s \la_t \right)_{s,t\in T},
\qquad
B= \left(|\lsp z_s, z_t \rsp|\right)_{s,t\in T}.
\end{equation}
Then the conditions (\ref{dl/dz}) and (\ref{dl/dla})
can be rewritten as
\begin{eqnarray}
\label{dl/dz_aa}
A \wt{z}_l  &=& \sum_{m=1}^n {\gamma}_{l m} \wt{z}_{m}
\quad {\rm for} \quad l = 1, \ldots, n\\
\label{dl/dla_bb}
B \wt{\la} &=& \bt  \wt{\la}.
\end{eqnarray}

Let $g= (g_{k l})_{k,l}$ be an $n \times n$ orthogonal  matrix  which
diagonalizes the hermitian $n \times n $  matrix
$\Gamma = (\gamma_{l m})_{l,m}$,
that  is, $g \,\Gamma g^* = D_\al$ is a diagonal
matrix with diagonal entries $\al_1, \ldots, \al_n$.

For $k = 1, \ldots, n$ set
$\wt{z}'_k = \sum_{l=1}^n g_{k l} \wt{z}_l$.
Then
$\wt{z}_m = \sum_{l=1}^n \overline{g}_{l m} \wt{z}'_l$,
for $m = 1, \ldots, n$.
We have $\lsp z'_s,  z'_t \rsp = \lsp z_s,  z_t \rsp$
for $1 \le s, t \le N$. Thus
the function $F$ and the matrices $A$ and $B$
do not change if we pass from variables induced by
the $\wt{z}_m$'s to the variables induced by
the $\wt{z}'_k$'s.
Similarly, the $\wt{z}'_k$'s satisfy the constraints
(\ref{constr_1}) and (\ref{constr_2}). Thus
the point $(z'_{s k}, \la_s)$ again gives a local
extremum of $F$, but with a new set of multipliers.

Expressing (\ref{dl/dz_aa}) in terms
of primed vectors $\wt{z}'_k$'s we get the $n$ eigenvalue
equations
\begin{equation}
A \wt{z}'_k = \al_k \wt{z}'_k
\quad \mbox{for\ }   k = 1, \ldots, n.
\label{dl/dz_bb}
\end{equation}
The last two conditions  mean that the  multipliers corresponding to
$(z'_{s k}, \la_s)$ are just
$\al_1, \ldots, \al_n, \bt$, with the off-diagonal
ones equal to 0.

Notice that if $\{f'{}^0_k\}$ is related to  $\{\wt{z}'_k\}$
by (\ref{subst}), then   $\{f'{}^0_k\}$ is an orthonormal basis
in $\spn[f_m^0]$; in particular the new square function $f'{}^0$
is equal to $f^0$.
Thus, without loss of generality, we can and will
work with these new ``primed''  vectors,
rather than with the original ones;
we will leave however the ``primes'' out,
for clarity of  notation. In other words, we will
assume  that
$(z_{s k}, \la_s)$ satisfies (\ref{constr_1}), (\ref{constr_2})
and (\ref{dl/dla_bb}), (\ref{dl/dz_bb}).

\medskip
We want to show that all $\al_k$'s are equal. To do so, we use
the well-known second order
conditions for a relative  maximum [H]:
the Hessian matrix $H$,
$$
H= \left( \begin{array}{ccc}
{\partial^2 L}/{\partial  Z_{p j} \partial  Z_{q k}} &\ &
{\partial^2 L}/{\partial  Z_{p j} \partial  \La_{q }} \\
 &  & \\
{\partial^2 L}/{\partial  \La_{p } \partial  Z_{q k}}&\ &
{\partial^2 L}/{\partial  \La_{p }\partial  \La_{q }}
\end{array}
\right),
$$
evaluated at the  point  $(z_{p j}, \la_p)$,
needs to be negative semi-definite on the tangent space
to the surface of constraints at that point.

%
We have
\begin{eqnarray}
  \label{dl/dzdz}
  \frac{\partial^2 L}{\partial  Z_{p j} \partial  Z_{q k}}
  &=&  \sgn \lsp z_p, z_q \rsp  \la_p \la_q \de_{jk}
  - \al_k \de_{pq} \de_{jk} \\
  \label{dl/dladla}
  \frac{\partial^2 L}{\partial  \La_{p} \partial  \La_{q }}
  &=& |\lsp z_p, z_q \rsp| - \bt \de_{pq}  \\
  \label{dl/dladz}
  \frac{\partial^2 L}{\partial  \La_{p } \partial  Z_{q k}}
  &=&  \sgn \lsp z_p, z_q \rsp  \la_q z_{p k} (1 + \de_{pq}).
\end{eqnarray}
So  $H$ is an ${N(n+1)}\times {N(n+1)}$ matrix of the form
\begin{equation}
H= \left( \begin{array}{cccc}
A - \al_1 I &  \ldots &0 & C_1^t\\
\vdots&  & \vdots& \vdots\\
0 & \ldots & A - \al_n I & C_n^t\\
C_1&  \ldots& C_n & B - \bt I
 \end{array} \right),
  \label{hessian}
\end{equation}
where $A$ and $B$ are defined in (\ref{a_b}),
$I$ is the identity matrix,
and    $C_k$ is the $N \times N$ matrix
$C_k = \left({\partial^2 L}
     /{\partial  \La_{p } \partial  Z_{q k}}\right)_{p,q=1}^N$
for $k=1, \ldots, n$.

The tangent space ${\cal T}$ to the surface of  constraints
described by (\ref{constr_1}) and  (\ref{constr_2})
consists of all vectors  $\wt{\wt{w}}  \in \Rn{N(n+1)}$,
\begin{equation}
\wt{\wt{w}}= \left( \begin{array}{c}
{w}_{pj}\\ \nu_p
\end{array}\right)
=  \left( \begin{array}{c}
\wt{w}_{1}\\\vdots\\\wt{w}_n\\ \wt{\nu}
\end{array}\right)
  \label{ww}
\end{equation}
orthogonal to all gradients $\nabla G_{0}$ and $\nabla G_{lm}$
for $1 \le l \le m \le n $,
evaluated at   $(z_{p j}, \la_p)$.
Then necessarily  $\lsp H \wt{\wt{w}}, \wt{\wt{w}}\rsp \le 0$
for all $\wt{\wt{w}} \in {\cal T}$.

>From (\ref{gradients}) it follows that
$\wt{\wt{w}}$ of the form (\ref{ww})
is in the tangent space $\cal T$ if and only if it
satisfies the following equations:
\begin{eqnarray}
   \hphantom{(3.22)\,}
 \lsp \nabla_{(Z_{pj},\La_p)} G_{lm} \wt{\wt{z}}, \wt{\wt{w}}\rsp
&=& \sum _{p\in T} (z_{pl} w_{pm} + z_{pm} w_{pl})
= \lsp \wt{z}_l, \wt{w}_m \rsp
   + \lsp \wt{z}_m, \wt{w}_l \rsp
=0\nonumber\\
\label{zz}
&&\qquad \mbox{ for\ } 1 \le l \le m \le n \\
  \lsp \nabla_{(Z_{pj},\La_p)} G_{0} \wt{\wt{z}}, \wt{\wt{w}}\rsp
&=& 2 \sum _{p\in T} \la_p \nu_p = 2 \lsp \wt{\la}, \wt{\nu} \rsp = 0.
\label{nu}
\end{eqnarray}

We will now show  that  for any
$1 \le l \ne m\le n$   and any $s \in T$
\begin{equation}
(\al_l - \al_m) z_{sl} z_{sm} = 0.
  \label{alam}
\end{equation}
This will follow from the negative-definiteness of
the  matrix $H$, by evaluating
$\lsp H \wt{\wt{w}}, \wt{\wt{w}}\rsp $
on suitable vectors $\wt{\wt{w}}$.

Consider the vector $\wt{\wt{w}^0}$ of the form (\ref{ww}) with
$\wt{w}_l = \wt{z}_m$, $\wt{w}_m = - \wt{z}_l$  and
$\wt{w}_k = 0$ otherwise,  and $\wt{\nu} \in \Rn{N}$
arbitrary satisfying (\ref{nu}). The orthogonality of the
vectors $\wt{z}_k$ ensured by (\ref{constr_1}) implies that
$\wt{\wt{w}^0} \in T$.

Let $\wt{H}$ denote the
$Nn \times Nn$ matrix which appears in the upper left corner of
(\ref{hessian}) and let $\wt{C}$ be the $N \times Nn$ matrix from
the bottom left corner.

A simple calculation using (\ref{dl/dz_bb}) shows that
$$
\lsp
   \left( \begin{array}{cc}
   \wt{H} & 0\\
    0 & 0
   \end{array} \right)
\wt{\wt{w}^0} , \wt{\wt{w}^0} \rsp =
\lsp
\left( \begin{array}{c}
(\al_m - \al_l)\wt{z}_l\\ (\al_m - \al_l) \wt{z}_m
\end{array}\right),
\left( \begin{array}{c}
\wt{z}_l\\ -\wt{z}_m
\end{array}\right)   \rsp =0.
$$
Thus
\begin{eqnarray}
\lsp H \wt{\wt{w}^0}, \wt{\wt{w}^0} \rsp
&=&
\lsp      \left( \begin{array}{cc}
   \wt{H} & \wt{C}^t\\
   \wt{C} & B  - \bt I
   \end{array} \right)
 \wt{\wt{w}^0}, \wt{\wt{w}^0} \rsp \nonumber\\
 &=&  \lsp (B - \bt I) \wt{\nu},  \wt{\nu}\rsp
 + \sum_{k=1}^n  \lsp C_k \wt{w}_k, \wt{\nu}\rsp
 + \sum_{k=1}^n  \lsp  \wt{w}_k, C_k^t \wt{\nu}\rsp
\nonumber\\
&=&
 \lsp (B- \bt I) \wt{\nu},  \wt{\nu}\rsp +
   2 \lsp  C_l \wt{z}_m - C_m \wt{z}_l, \wt{\nu}\rsp \le 0.
  \label{definit}
\end{eqnarray}


Replacing  $\wt{\nu}$ by $\ep \wt{\nu}$
and taking $\ep \to 0$ we get that the first
term in the final sum (which is  of the
second order in $\ep$)
can be disregarded from the estimate.
In the  inequality obtained this way
$\wt{\nu}$ can be replaced by
$-\wt{\nu}$, hence
$ \lsp  C_l \wt{z}_m - C_m \wt{z}_l, \wt{\nu}\rsp =0$.
Since this equality holds for an
arbitrary  $\wt{\nu}$
orthogonal to $\wt{\la}$,
we conclude that there is a constant $\gamma $
such that  $ C_l \wt{z}_m - C_m \wt{z}_l = \gamma \wt{\la}$.
Equivalently, the definition of $C_l$ and
(\ref{dl/dladz}) yield that for $p \in T$,
\begin{eqnarray}
\gamma \la_p &=&
 z_{p l}\sum_{q\in T} \sgn \lsp z_p, z_q \rsp  \la_q z_{q m} (1 + \de_{pq})
- z_{p m}\sum_{q\in T} \sgn \lsp z_p, z_q \rsp  \la_q z_{q l} (1 + \de_{pq})
\nonumber\\
&=&  z_{p l}\sum_{q\in T} \sgn \lsp z_p, z_q \rsp  \la_q z_{q m}
- z_{p m}\sum_{q\in T} \sgn \lsp z_p, z_q \rsp  \la_q z_{q l}.
\label{ccc}
\end{eqnarray}
Observe that  by  (\ref{a_b}) and (\ref{dl/dz_bb}) we have,
for $p\in T$,
$$
\sum_{q\in T} \sgn \lsp z_p, z_q \rsp \la_p \la_q z_{q m}
= (A \wt{z}_m)_p = \al_m z_{pm}
$$
and an analogous equality holds for the second term in
(\ref{ccc}). This implies  that
\begin{equation}
\gamma \la_p^2 = ( \al_m - \al_l) z_{pm}  z_{pl}
\quad {\rm for} \quad p \in T.
  \label{ggg}
\end{equation}
Summation over $p$ yields by
(\ref{constr_1}) and (\ref{constr_2}) that
$$
\gamma = ( \al_m - \al_l)  \lsp\wt{z}_m, \wt{z}_l\rsp =0.
$$
Hence (\ref{alam}) holds.
It obviously follows from Lemma~\ref{chain} that
\begin{equation}
 \al_m = \al_l =: \al \qquad {\rm for \ all\ }  1 \le l, m \le n.
  \label{all_alam}
\end{equation}
Expressing (\ref{dl/dz_bb}) coordinatewise
we have
$$
\sum _{t\in T} \sgn \lsp z_s, z_t \rsp  z_{t l} \la_s \la_t
- \al z_{s l} =0 \quad {\rm for\ } l=1, \ldots, n,\ s \in T.
$$
Multiplying by $ z_{s l}$, summing up over $l$ and using
(\ref{dl/dla_bb}),  we find, for $s\in T$,
$$
0 = \sum _{t\in T} | \lsp z_s, z_t \rsp| \la_s \la_t
- \al \sum _{l=1}^n z_{s l}^2
= \bt \la_s^2  - \al \sum _{l=1}^n z_{s l}^2.
$$

In terms of $\mu_s$ and $f^0$ this means that
$\bt \mu_s = \al f^0(s)\mu_s$,
\ie\  $f^0(s) = \bt / \al$ is constant for
all $s \in T$ with $\mu_s \ne 0$. If $\mu_s =0$, then
$f^0(s)=0$, as mentioned already.
Since
$\sum_{s\in T} |f^0(s)|^2 \mu_s =n = (\bt / \al)^2$,
we conclude that $\bt / \al = \sqrt n$,
completing the proof in the  case (a).

\medskip
\noindent{\bf (b)}
We now consider the infinite case $T = \NN$ for $\Kn{}= \Rn{}$.
Assume that the function $F$ given by (\ref{obj_funct}) attains a
relative maximum subject to the constraints (\ref{constr_1})
and (\ref{constr_2}) at the point $ (z_{s j}, \la_s)$,
where $z_{sj} = f^0_j(s) \sqrt{\mu^0_s}$
and $\la_s = \sqrt{\mu^0_s}$ for $s \in T$,
$ j = 1, \ldots, n$.
Note that $\wt{z}_j = ( z_{s j})_{s \in T} \in l_2$,
for $ j = 1, \ldots, n$, and
$\wt{\la} = (\la_s)_{s \in T}  \in l_2$.
For $N \in \NN$ and an arbitrary vector $\wt{z} \in l_2$
set $\wt{z}^N = R_N \wt{z}\in l_2^N$.

Fix $N \in \NN$ sufficiently large so that
$\wt{z}_1^N, \ldots, \wt{z}_n^N$ are linearly
independent. In  (\ref{obj_funct})--(\ref{constr_2})
fix the variables for  $s >N$ by putting
$$
Z_{sj} = z_{sj}, \  \Lambda_s = \lambda_s \quad
\mbox{for\ } s > N, j=1, \ldots, n.
$$
Relative to the new constraints, $F$ as a function in the
variables $(Z_{s j}, \La_s)$,
with $s= 1,\ldots, N$, $j=1,\ldots, n$,
attains a relative maximum at $ (z_{s j}, \la_s)$, with
$s= 1,\ldots, N$, $j=1,\ldots, n$.
The first order Lagrange multiplier conditions
(\ref{dl/dz}) and (\ref{dl/dla}) now take the form
\begin{eqnarray}
\label{inf_dl/dz}
\sum _{t\in T} \sgn \lsp z_s, z_t \rsp  z_{t l} \la_s \la_t
&-& \sum_{m=1}^n {\gamma}_{l m} z_{s m} =0 \nonumber\\
&&  \mbox{for\ }  s= 1,\ldots, N, l=1, \ldots, n \\
\label{inf_dl/dla}
\sum _{t\in T} |\lsp z_s, z_t \rsp| \la_t
&-& \bt  \la_s =0 \quad \mbox{for\ } s= 1,\ldots, N.
\end{eqnarray}

For a fixed $ s= 1,\ldots, N$, the first sum in
(\ref{inf_dl/dz})  is independent of $N$.
Since $\wt{z}_1^N, \ldots, \wt{z}_n^N$ are linearly
independent, this uniquely  determines the matrix
$\Gamma = (\gamma_{l m})_{l, m =1}^n$, which is then
independent of $N$. Thus
(\ref{inf_dl/dz}) and (\ref{inf_dl/dla}) hold
for all $s \in \NN$.
Again, we diagonalize the $n \times n$ matrix $\Gamma$
and we introduce $\wt{z}_1', \ldots, \wt{z}_n'$
satisfying the eigenvalue equations
$A \wt{z}'_k = \al_k \wt{z}'_k$ for
$ k = 1, \ldots, n$. Moreover,
$B \wt{\la} = \bt  \wt{\la}$.
Note that $A$ and $B$, formally given by
(\ref{a_b}), are now infinite Hilbert--Schmidt
matrices. Again, in what follows, we leave the
``primes'' out and write simply $\wt{z}_k$.

Now let $\cal T$  be the space of (infinite) vectors
$\wt{\wt{w}}$ in the direct sum $ \bigoplus l_2$
of $n+1$ copies of $l_2$, which are of the form
(\ref{ww}) and satisfy (\ref{zz}) and (\ref{nu})
(with $T = \NN$).
Let $H$ be the (infinite) matrix of the form
(\ref{hessian}), with $A$, $B$ and $C_k$ being
infinite as well. To conclude the same proof as in part (a),
it suffices to show that
$\lsp H \wt{\wt{w}}, \wt{\wt{w}}\rsp \le 0$ for
all $\wt{\wt{w}} \in {\cal T}$.

To this end, denote by $H^N$ and $A^N$, $B^N$, $C_k^N$
the restricted matrices of order $N (n+1)\times N (n+1)$
and $N \times N$ respectively.

The constraints for the restricted problem in the variables
$(Z_{s j},\La_s)$,
with $s= 1,\ldots, N$, $j=1,\ldots, n$, are still of the form
\begin{eqnarray*}
  G_{l m}^N(Z_{s j}, \La_s) &=&
  \sum _{s=1}^N Z_{s l} \overline{Z}_{s m} - d_{l m} = 0 \\
  G_0^N (Z_{s j}, \La_s) &=&
        \sum_{s=1}^N \La _s ^2 - d =0,
\end{eqnarray*}
for some $d, d_{l m} \in \Rn{}$. This implies that the
corresponding tangent space ${\cal T}^N \subset \Rn{N (n+1)}$
of vectors $\wt{\wt{w}}_N$ of the form (\ref{ww})
is defined by the equations
\begin{eqnarray}
\label{Nzz}
  \sum _{p=1}^N (z_{pl} w_{pm} + z_{pm} w_{pl})
&=&0 \qquad \mbox{for\ } 1 \le l \le m \le n \\
 \sum _{p=1}^N \la_p \nu_p &=&  0.
\label{Nnu}
\end{eqnarray}
Hence, in general, the projection
$\wt{\wt{w}}^N = R_{N(n+1)}\wt{\wt{w}}$
of $\wt{\wt{w}}$ onto $\Rn{N(n+1)}$ is not in
${\cal T}^N$, since
$\lsp \wt{z}_l, \wt{w}_m \rsp  + \lsp \wt{z}_m, \wt{w}_l \rsp =0$
for  $\le l \le m \le n $ does not imply
$\lsp \wt{z}_l^N, \wt{w}_m^N \rsp  +
\lsp \wt{z}_m^N, \wt{w}_l ^N\rsp =0$
for  $\le l \le m \le n $.
However, since  the limit of
(\ref{Nzz}) and (\ref{Nnu}),  as $N\to \infty$,
coincides with (\ref{zz}) and (\ref{nu})
(for $T = \NN$),
it is clear that for any $\wt{\wt{w}} \in {\cal T}$,
there is a sequence $(\wt{\wt{w}}_N)_{N=1}^\infty$,
with
$$
\wt{\wt{w}}_N
= \left( \begin{array}{c}
(\wt{w}_N)_{l}\\ \wt{\nu}_N
\end{array}\right)
\ \in {\cal T}^N,
$$
such  that  $\wt{\wt{w}}_N \to \wt{\wt{w}}$ in the
$\bigoplus l_2$-norm.

Since $\lsp H^N \wt{\wt{w}}_N, \wt{\wt{w}}_N\rsp \le 0$
for all $N \in \NN$, it suffices to show that
$$
\lim_{N\to \infty} \lsp H^N \wt{\wt{w}}_N, \wt{\wt{w}}_N\rsp
= \lsp H \wt{\wt{w}}, \wt{\wt{w}}\rsp.
$$
This is shown term by term. A typical case is
$$
\lim_{N\to \infty}
\lsp (A^N - \al_l I)(\wt{w}_N)_l, (\wt{w}_N)_l\rsp
= \lsp A \wt{w}_l - \al_l  \wt{w}_l, \wt{w}_l \rsp,
$$
which reduces to
\begin{equation}
\lsp (\wt{w}_N)_l, (\wt{w}_N)_l\rsp \to
\lsp  \wt{w}_l, \wt{w}_l \rsp
\mbox{\ and \ }
\lsp A^N(\wt{w}_N)_l, (\wt{w}_N)_l\rsp \to
\lsp A \wt{w}_l, \wt{w}_l \rsp.
\label{lim}
\end{equation}
But (\ref{lim}) follows from
$\lim_{N\to \infty} (\wt{w}_N)_l = \wt{w}_l$
in the $l_2$-norm and the fact that matrices $A^N$ converge
to $A$ in the Hilbert--Schmidt norm.

As before, we find that $\al_1 = \ldots = \al_n =:\al$
and $\bt \mu_s = \al f^0(s)\mu_s$. We then complete the proof
as in case (a).

\medskip
\noindent{\bf (c)}
Finally, we indicate the necessary changes in the proof of the
complex case, $\Kn{}= \Cn{}$.
We assume for simplicity that $T = \{1, \ldots, N\}$.
The function $F$, as defined by (\ref{obj_funct}), is now
a function of the complex variables
$Z_{s j}= X_{s j} + i Y_{s j}$ and the real variables
$\La_s$. We consider $F$ as a function of real variables
$(X_{s j}, Y_{s j}, \La_s)$.
There are now $ n^2$ real constraints for
$1 \le l \le m \le n$,
\begin{eqnarray}
G_{l m}^{(1)}(X_{s j},Y_{s j}, \La_s)&:=& \RRe G_{l m}(Z_{s j}, \La_s)=0
\mbox{\quad for\ }l < m \nonumber\\
\label{constr_1_c}
G_{l m}^{(2)}(X_{s j},Y_{s j}, \La_s)&:=&
         \IIm G_{l m}(Z_{s j}, \La_s)=0 \mbox{\quad for\ }l < m \\
G_{l l}(X_{s j},Y_{s j}, \La_s)&:=& G_{l l}(Z_{s j}, \La_s)=0
\mbox{\quad for\ }l = m,\nonumber
\end{eqnarray}
as well as  (\ref{constr_2}).

Consider the Lagrange function $L$ defined by
$$
2 L := F -
\sum_{l < m} (\wt{\gamma}_{l m}^{(1)} G_{l m}^{(1)} -
\wt{\gamma}_{l m}^{(2)} G_{l m}^{(2)})
- \sum_{l } \wt{\gamma}_{l l} G_{l l} -
\bt G_0.
$$
If $F$ attains the extremum subject to conditions
(\ref{constr_2}) and (\ref{constr_1_c}) at
$(z_{s j}= x_{s j}+i y_{s j}, \la_s )$, then a calculation
shows that the first order conditions can be written
in the following complex form
\begin{equation}
\frac{\partial L}{\partial  X_{s l}} +
i \frac{\partial L}{\partial  Y_{s l}} =
\sum _{t\in T} \sgn \lsp z_s, z_t \rsp  z_{t l} \la_s \la_t
- \sum_{m=1}^n {\gamma}_{l m} z_{s m} =0,
\label{c_dl/dz}
\end{equation}
for $s \in T$, $l=1, \ldots, n$.
Here $\sgn w = w / |w|$ for $w \in \Cn{}$, $w \ne 0$
and  $\sgn 0 =0$.
Moreover,
$$
\gamma_{l m} =   \frac{1}{2}  \left\{
\begin{array}{ll}
  \wt{\gamma}_{l m}^{(1)} - i \wt{\gamma}_{l m}^{(2)}
                & \mbox{if  $l< m$}\\
   \wt{\gamma}_{ m l }^{(1)} + i  \wt{\gamma}_{ m l }^{(2)}
                & \mbox{if  $m< l$}\\
                      2 \wt{\gamma}_{ l l } & \mbox{if  $m= l$}
\end{array}
 \right.
$$
define an $n \times n$ hermitian complex matrix $\Gamma$.
Thus we can again diagonalize $\Gamma$ and rewrite
(\ref{c_dl/dz})  and
(\ref{dl/dla}) as
\begin{equation}
A \wt{z}_k = \al_k \wt{z}_k
\quad \mbox{for\ }   k = 1, \ldots, n
\quad  \mbox{and\ }
B \wt{\la} = \bt  \wt{\la},
\label{1_order_c}
\end{equation}
where $A$ and $B$ are formally defined as in (\ref{a_b}).

The tangent space $\cal T$ to the surface of constraints
(\ref{constr_1_c}) and (\ref{constr_2}) now consists of
vectors $\wt{\wt{w}} \in \Cn{N(n+1)}$, whose complex form
is formally described by (\ref{ww}) and whose real form is
\begin{equation}
\wt{\wt{w}} =
\left( \begin{array}{c}
\wt{u}_{1}\\\wt{v}_{1}\\ \vdots\\ \wt{u}_n\\ \wt{v}_n\\ \wt{\nu}
\end{array}\right)
\in \Rn{N(2n+1)} 
\label{ww_c}
\end{equation}
where $\wt{w}_{l} = \wt{u}_{l} + \wt{v}_{l}$ for
$l=1, \ldots, n$.
The equations defining $\cal T$ can be written in
the following (complex) form
\begin{eqnarray}
\label{zz_c}
\sum _{p\in T} (z_{pl} \overline{w}_{pm} + w_{pl} \overline{z}_{pm})
&=& \lsp \wt{z}_l, \wt{w}_m \rsp
   + \lsp \wt{w}_l, \wt{z}_m  \rsp  =0 \nonumber\\
&& \mbox{ for\ } 1 \le l \le m \le n,\\
\label{nu_c}
\lsp \wt{\la}, \wt{\nu} \rsp &=& 0.
\end{eqnarray}

The Hessian matrix $H$ in the real form has now size
$N (2n +1) \times N (2n +1)$. In particular, the
matrix $C$ in (\ref{definit}) consists of $2n$ real
matrices
$$
C_l^{(1)} = \left(\frac{\partial^2 L}
     {\partial  \La_{p } \partial  X_{q l}}\right)_{p,q=1}^N,
\quad
C_l^{(2)} = \left(\frac{\partial^2 L}
     {\partial  \La_{p } \partial  Y_{q l}}\right)_{p,q=1}^N
\quad\mbox{for\ } l=1, \ldots, n
$$
of size $N \times N$, evaluated at $(x_{sj}, y_{sj}, \la_s)$.
The condition  $\lsp H \wt{\wt{w}}, \wt{\wt{w}}\rsp \le 0$
translates into
\begin{equation}
\sum_{l=1}^n \lsp C_l^{(1)} \wt{u}_l +  C_l^{(2)} \wt{v}_l,
     \wt{\nu}\rsp =0,
\label{neg_def_c}
\end{equation}
for all $\wt{\nu}$ satisfying (\ref{nu_c}).
Thus
$\sum_{l=1}^n ( C_l^{(1)} \wt{u}_l +  C_l^{(2)} \wt{v}_l)$
is a multiple of $\wt{\la}$, for all
$\wt{\wt{w}}\in {\cal T}$ of the (complex) form  (\ref{ww})
satisfying  (\ref{zz_c}) and (\ref{nu_c}).

For $1 \le l \ne m \le n$ we pick two different types of vectors
in $ {\cal T}$. In the complex form  (\ref{ww}), the first
vector $\wt{\wt{w}}$ looks as before, that is,
$\wt{w}_l = \wt{z}_m$, $\wt{w}_m = - \wt{z}_l$  and
$\wt{w}_k = 0$ otherwise,  and $\wt{\nu} \in \Rn{N}$
arbitrary satisfying (\ref{nu_c}).
The second type,  $\wt{\wt{w}}'$, is defined similarly by setting
$\wt{w}_l' = i\wt{z}_m$, $\wt{w}_m' = i \wt{z}_l$  and
$\wt{w}_k = 0$ otherwise.
The real form of these vectors is the following, writing
the non-zero terms only,
$$
\wt{\wt{w}} =
\left( \begin{array}{c}
\wt{x}_{m}\\ \wt{y}_{m}\\- \wt{x}_l\\ - \wt{y}_l\\ \wt{\nu}
\end{array}\right)
\qquad \mbox{and}\qquad
\wt{\wt{w}}' =
\left( \begin{array}{c}
-\wt{y}_{m}\\ \wt{x}_{m}\\ - \wt{y}_l\\ \wt{x}_l\\ \wt{\nu}
\end{array}\right).
$$
Both vectors satisfy  (\ref{zz_c}) and (\ref{nu_c}).
Calculating (\ref{neg_def_c}) for  $\wt{\wt{w}}$ and $\wt{\wt{w}}'$
and using the eigenvalue equations (\ref{1_order_c}) we find,
in an analogous  way as we obtained (\ref{ggg})
in case (a), that there are $\gamma_1$
and $\gamma_2$ such that
\begin{eqnarray*}
\la_p \Bigl( (C,0)\wt{\wt{w}}\Bigr)_p &=&
\RRe ( \al_m - \al_l) \overline{z}_{pm}  z_{pl}
=  \gamma_1 \la_p^2\\
\la_p \Bigl( (C,0)\wt{\wt{w}}'\Bigr)_p &=&
\IIm ( \al_m - \al_l) \overline{z}_{pm}  z_{pl}
=  \gamma_2 \la_p^2
\end{eqnarray*}
for all $ p \in T$.
Summing over $p$ and using (\ref{constr_1_c}) and (\ref{constr_2})
we infer that $\gamma_1= \gamma_2=0$, thus
$$
( \al_m - \al_l)  \overline{z}_{pm}  z_{pl}=0.
$$
Just as in case (a), the last equality implies
$\al_1 = \ldots = \al_n =:\al$ and
$\bt \mu_s = \al f^0(s) \mu_s$ for all $s \in T$,
which completes the proof of (c).
\qed

\section{The estimate for the projection constant}

We start by a simple but useful lemma of Sidelnikov [Si] and Goethals
and Seidel [GS.1].  It gives a lower bound for expressions related to
those appearing in the definition (\ref{sups}) of $\phi$. Since the
bound is essential for our estimate, we include its proof.

\begin{lemma}
  \label{curvature}
  Let $T = \{1, \ldots, N\}$, or $T = \NN$ and let $(\mu_s)_{s \in T}$
  be a probability measure on $T$.  Let $(z_s)_{s \in T} \in \Kn{n}$
  with $\|z_s\|_2 =1$.  Let $\omega$ be the normalized
  rotation-invariant measure on $S^{n-1}= S^{n-1}(\Kn{})$.  Then for
  every even natural number $k \in 2\NN$,
  \begin{equation}
\sum_{s, t \in T} |\lsp z_s, z_t\rsp|^{k}\mu_s \mu_t
\ge \int_{S^{n-1}}\, \int_{S^{n-1}}\,
     |\lsp z, w\rsp|^{k}\, d\omega (z) d\omega (w).
    \label{curv}
  \end{equation}
\end{lemma}

(In the complex case, express the integrand in the real variables
and integrate over   $S^{n-1}(\Cn{})= S^{2n -1}(\Rn{}) $.)

\medskip
\noindent\proof
Let $n \in \NN$ and $k=2m \in 2\NN$. For $z \in \Kn{n}$, let
$z^{\otimes j}= z \otimes \ldots  \otimes z \in \Kn{n^j}$
denote the $j$-fold tensor product of $z$ with itself,
for $j = 1, 2,\ldots$.
Scalar products in $\Kn{n^j}$ will be denoted by
$\lsp \cdot,  \cdot\rsp _j$, and for $j=1$ just by
$\lsp \cdot,  \cdot\rsp $. Then
for any $z, w \in \Kn{n}$ and $j = 1, 2, \ldots$
we have
$$
\lsp z^{\otimes j}, w^{\otimes j}\rsp_j = \lsp z, w\rsp ^j,
$$
and
$$
\lsp z^{\otimes m} \otimes \overline{z}^{\otimes m},
   w^{\otimes m} \otimes \overline{w}^{\otimes m} \rsp_k =
    \lsp z, w\rsp ^m \, \lsp \overline{z},\overline{ w}\rsp ^m
= |  \lsp z, w\rsp |^ k.
$$
Consider
$$
\xi:= \sum_{s \in T} (z_s^{\otimes m} \otimes
\overline{z}_s^{\otimes m})
\mu_s -  \int_{S^{n-1}}\,
(z^{\otimes m} \otimes \overline{z}^{\otimes m})  \, d\omega (z)
\in \Kn{{n^k}}.
$$
By the rotation invariance of $\omega$, integrals of the form
$ \int_{S^{n-1}}\, |\lsp e, w\rsp|^{k}\, d\omega (w)$
do not depend on $e \in S^{n-1}$. This allows to evaluate
$\lsp \xi, \xi \rsp_k$ as follows:
\begin{eqnarray*}
0 \le \lsp \xi, \xi \rsp_k
&=&  \sum_{s, t \in T} |\lsp z_s, z_t\rsp|^{k}\mu_s \mu_t
+ \int_{S^{n-1}}\, \int_{S^{n-1}}\,
     |\lsp z, w\rsp|^{k}\, d\omega (z) d\omega (w)\\
&& \qquad\qquad -2 \sum_{s \in T}\mu_s  \int_{S^{n-1}}\,
     |\lsp z_s, w\rsp|^{k}\, d\omega (w)\\
&=&  \sum_{s, t \in T} |\lsp z_s, z_t\rsp|^{k}\mu_s \mu_t
- \int_{S^{n-1}}\, \int_{S^{n-1}}\,
     |\lsp z, w\rsp|^{k}\, d\omega (z) d\omega (w),
\end{eqnarray*}
which proves the lemma.
\qed

Now we are ready for the proof  of Theorem~\ref{main}.

\smallskip
\noindent{\bf Proof of Theorem~\ref{main}\ } {\bf (a)}
Let $n \in \NN$ and let
$G(n)$ denote the right hand side of
(\ref{ineq_main}). We have to show that for any $n$-dimensional
space $E$ we have $\la (E) \le G(n)$.
By Proposition~\ref{lewis},
$\la (E) \le \phi (n, T)$, where
$T= \{1, \ldots, N\}$, if $E \subset l_\infty^N$,
and $T = \NN$, if $E \subset l_\infty$.
By Lemma~\ref{phi_finite} it suffices
to show that
\begin{equation}
\phi (n, T) \le G(n) \qquad \mbox{for\ }
T= \{1, \ldots, N\}.
\label{phi_g}
\end{equation}

Given $n$ and $T$ we may assume that
$\phi (n_1, T_1) < \phi (n, T) $ for all
$n_1 < n$ and all $T_1 \subset T$; otherwise
the proof which follows would be applied to the minimal
$n_1$ with $\phi (n_1, T) = \phi (n, T) $, to show
that
$\phi (n, T) \le G(n_1) < G(n)$.

The double supremum in the definition (\ref{sups}) of $\phi(n,T)$
is attained for some probablity measure
$\mu = (\mu_s)_{s\in T}$ on $T$ and some
orthonormal system $f_j= (f_j(s))_{s\in T} \in l_2(T, \mu)$,
$j=1, \ldots, n$. By Proposition~\ref{sq_funct},
the square function $f:= (\sum_{j=1}^n |f_j|^2)^{1/2}$
equals $\sqrt n$ for all $s$ where $\mu_s \ne 0$,
and equals to $0$ otherwise. The span of the $f_j$'s is
therefore supported by $S:= \mbox{supp} \mu \subset T$.
For $s \in S$, let
$z_s:= n^{-1/2}(f_j(s))_{j=1}^n \in l_2^n$.
Then $\|z_s\|_2 =1$ and
\begin{equation}
\phi(n, T) =  n \sum_{s, t \in S} |\lsp z_s, z_t \rsp | \mu_s \mu_t.
\label{max}
\end{equation}

Define $\al$ and $\bt$ by
$$
\al; = \left\{
         \begin{array}{ll}
1/ \sqrt{n+2} &\ \Kn{} = \Rn{}\\
1/ \sqrt{n+1} &\ \Kn{} = \Cn{}
        \end{array}
\right.,
\qquad
\bt: = \left\{
         \begin{array}{ll}
3/ (n+2) &\ \Kn{} = \Rn{}\\
2/(n+1) &\ \Kn{} = \Cn{}
        \end{array}
\right..
$$
Then for $u \in [-1, 1]$ we have
$$
(|u| - \al)^2 = \Bigl( (u^2 - \al^2)/ (|u| + \al) \Bigr)^2
\ge  (u^2 - \al^2)^2 / (1+\al)^2.
$$
This implies
\begin{equation}
|u| \le \gamma_0 +\gamma_2 u^2 - \gamma_4 u^4
\qquad \mbox{for\ } u \in [-1, 1],
\label{gamma_ineq}
\end{equation}
where
\begin{equation}
\gamma_0 = \frac{\al}{2} - \frac{\al^3}{2(1+\al)^2},
\ \gamma_2 = \frac{1}{2\al} + \frac{\al}{(1+\al)^2},
\ \gamma_4 = \frac{1}{2\al (1+\al)^2}
\label{gammas}
\end{equation}
are non-negative. Equality in (\ref{gamma_ineq}) occurs
for $ u \in [-1, 1]$ if and only if $|u|$ equals
to 1 or $\al$. (The right hand side of (\ref{gamma_ineq})
touches $|u|$ at $\pm \al$ and intersects
$|u|$ at $\pm 1$.)

Using (\ref{gamma_ineq}) and (\ref{curv}) we can estimate
(\ref{max}).
\begin{eqnarray}
\phi(n,T) &\le&  n  \sum_{s,t \in T}
\Bigl( \gamma_0 + \gamma_2 |\lsp z_s, z_t \rsp|^2
   - \gamma_4 |\lsp z_s, z_t \rsp|^4 \Bigr)\mu_s \mu_t \nonumber\\
&\le & n \Bigl( \gamma_0 + \gamma_2/n - \gamma_4
\int_{S^{n-1}}\, \int_{S^{n-1}}\,
     |\lsp z, w\rsp|^{4}\, d\omega (z) d\omega (w)\Bigr)\nonumber\\
& = & n \gamma_0 + \gamma_2 - \gamma_4 \bt = G(n).
\label{fin_est}
\end{eqnarray}

Here we used the orthonormality of the $f_j$'s to evaluate
the double sum
$$
\sum_{s,t \in T} |\lsp z_s, z_t \rsp|^2 \mu_s \mu_t = 1/n
$$
and the fact that for any $e \in S^{n-1}$,
$$
I:=  \int_{S^{n-1}}\,
     |\lsp e, w\rsp|^{4}\,  d\omega (w)= \bt / n,
$$
since \eg\ in the real case,
$$
I =  \int_{-1}^1 t^4 (1-t^2)^{(n-3)/2}\, dt /
   \int_{-1}^1  (1-t^2)^{(n-3)/2}\, dt = 3/(n (n+2));
$$
in the complex case the calculation yields
$I = 2/ (n (n+1))$.

The last equality in (\ref{fin_est}) is established by
a direct calculation using (\ref{gammas}).

\medskip
{\bf (b)} and {\bf (c)}
We now assume that $E$ is an $n$-dimensional space attaining
the extremal bound, $\la(E) = G(n)$.
By Proposition~\ref{lewis}, there is
$T = \{1, \ldots, N\}$ or $T = \NN$, a probability measure
$\mu = (\mu_s)_{s \in T}$ on $T$ and an orthonormal  basis
$(f_j)_{j=1}^n$ in $E_{2,\mu}$ such that
\begin{equation}
\la (E) \le \sum_{s, t\in T} | \sum_{j=1}^n f_j(s)
                   \overline{f_j (t)}| \mu_s \mu_t.
\label{extr_case}
\end{equation}

For all $n_1 < n$ and all $T_1 \subset T$ we have
$\phi (n_1, T_1) < \phi (n, T)$;
otherwise, for some  $n_1 < n$ and some
$T_1 \subset T$ we would have, by part (a),
$$
G(n) = \la(E) \le \phi (n, T) \le
\phi (n_1, T_1) \le G(n_1) < G(n).
$$
Therefore, by Proposition~\ref{sq_funct}, on the support
$S \subset T$ of $\mu$, the square function
$f= (\sum_{j=1}^n |f_j|^2)^{1/2}$ is equal to $\sqrt n$.
For $s \in S$ consider again the short vectors
$z_s = (f_j(s))_{j=1}^n / \sqrt{n}$. Hence
$\|z_s\|_2 =1$ for $s \in S$.
We may and will further assume that $S$ is minimal
in the sense that for $s \ne t$ we have
$z_s \ne \theta z_t$ with $|\theta |=1$.
Otherwise, we could replace the short vectors $z_s$
and $z_t$ by one vector $z_s$, assigning to it
the measure $\mu_s + \mu_t$; the orthogonality
and the normalization
of the corresponding long vectors  and the double sum
in (\ref{extr_case}) would remain unchanged.
Let $N:= |S|$. We have to show that $N$ is finite,
and, in fact, bounded by $N(n)$ as defined in (\ref{Nn}).

By  (\ref{extr_case}), (\ref{gamma_ineq}) and (\ref{fin_est})
we have
\begin{eqnarray*}
\la (E) &=& n \sum_{s, t\in S} |\lsp z_s, z_t \rsp| \mu_s \mu_t\\
&\le& n \sum_{s, t\in S}
\Bigl(\gamma_0 + \gamma_2 |\lsp z_s, z_t \rsp|^2
   - \gamma_4 |\lsp z_s, z_t \rsp|^4 \Bigr)\mu_s \mu_t \\
&\le& G(n).
\end{eqnarray*}
Thus, the assumption  $\la (E)= G(n)$ implies the equality
of all  terms.  The equality in the first inequality
requires that $|\lsp z_s, z_t\rsp| = \al$ or 1 for all
$s, t \in S$ (note that $\mu_s \ne 0$ for $s \in S$).
For $s \ne t$, $z_s \ne \theta z_t$, hence
$|\lsp z_s, z_t\rsp| = \al$.
Recall that $\al = 1/\sqrt{n+2}$ in the real case, and
$\al =  1/\sqrt{n+1}$ in the complex case.
We thus proved that the vectors
$(z_s)_{s \in S} \subset S^{n-1}(\Kn{})$
are equiangular. Since in $\Kn{n}$ there are
at most $N(n)$ equiangular vectors, it follows that
$N= |S| \le N(n)$.
Using the Cauchy--Schwartz inequality, we get another chain
of inequalities which become   equalities,
\begin{eqnarray*}
G(n) &=& n\Bigl( \sum_{s, t \in S} \mu_s \mu_t \al +
       \sum_{s \in S} \mu_s^2 (1 -\al)\Bigr)\\
&\ge& n \al + (n/N) (1 - \al) \ge
   n \al + (n/N(n)) (1 - \al) = G(n),
\end{eqnarray*}
where the last equality follows by  a direct calculation,
inserting the value of $\al$. The equality implies, in particular,
that $N = N(n)$. Also, all values of $\mu_s$ have to be equal
($\mu_s = 1/N(n)$). Since the vectors $f_j$ are all supported
by $S$, it follows that $E$ is isometric to a subspace
of $l_\infty^{N(n)}$.
The orthogonal projection, given by  the matrix
$(n / N(n)) (\lsp z_s, z_t \rsp )_{s,t}$, is a minimal
projection.

Conversely, if in $\Kn{n}$ there exist $N(n)$ equiangular vectors
$(z_s)$, we may construct $E = \spn [f_1, \ldots, f_n]
\subset l_\infty^{N(n)}$  by letting $f_j(s) = \sqrt{n}\, z_{sj}$
($j=1, \ldots, n$, $s= 1, \ldots, N(n)$).
By [K], the projection constant of $E$ is equal to $G(n)$,
and the $f_j$'s are orthonormal with respect to the
equidistributed  probability measure $\mu$ on $\{1, \ldots, N(n)\}$.
Moreover,  $P$, given by the matrix
$(n/N(n)) (\lsp z_s, z_t \rsp )_{s, t=1}^{N(n)}$,
and acting as an operator from  $l_\infty^{N(n)}$
to  $l_\infty^{N(n)}$,
is a minimal (and orthogonal)
projection onto $E$  with norm $G(n)$.

Either way, the norm of a vector $ \sum_{j=1}^n \al_j f_j $
in $l_\infty^{N(n)}$ is given by
$$
\|\sum_{j=1}^n \al_j f_j\|_\infty =
\sup_{1 \le s \le N(n)} \sqrt{n}\, |\lsp \al, z_s \rsp|.
$$
Thus, given $N(n)$ equiangular vectors $(z_s)$ in $\Kn{n}$,
we get an $n$-dimensional normed space with the maximal projection
constant by setting
\begin{equation}
\|(\al_j)_{j=1}^n\| := \sup_{1 \le s \le N(n)} |\lsp \al, z_s \rsp|.
\label{norm}
\end{equation}

\medskip

In the real case, $N(n) = n(n+1)/2$ equiangular vectors exist
in $\Rn{n}$ for $n = 2, 3, 7, 23$ and  these systems are unique
up to orthogonal transformations.
Hence the real spaces with  projection constant $G(n)$
are unique up to isometry if  $n = 2, 3, 7, 23$.
For $n=2$, the uniqueness (up to orthogonal transformations)
of the three vectors at angle $2\pi/3$ each, is trivial.
For $n=3$ one considers the $6 \times 6$ Gram matrix
$(\lsp z_s, z_t\rsp)$, with
$|\lsp z_s, z_t\rsp|= 1 / \sqrt 5$ for $s \ne t$.
It is easy to see that up to permutations and multiplications
of the $z_s$'s by $-1$, the sign pattern is uniquely
determined. The standard paper on the subject is
Lemmens, Seidel [LS]; for the uniqueness for $n=7, 23$
we refer to Goethals, Seidel [GS.2] and Seidel [S].
For $n=2$, (\ref{norm}) yields the norm with the (regular) hexagonal
unit ball; for $n=3$, the extremal ball defined via
(\ref{norm}) is the (regular) dodecahedron, since the 6
equiangular vectors in $\Rn{3}$ are the diagonals
of the icosahedron.

In the complex case, $N(n)= n^2$ equiangular vectors
exist in  $\Cn{n}$ at least for $n=2, 3$. For $n=2$, the system
and the extremal space are again unique up to isometry.
For $n=3$, the system of $9$ vectors in $\Cn{3}$
is connected to the Hessian polyhedron, \cf\
Coxeter [C].
\qed

\rem
Part (a) of the previous proof also shows that
Theorem~\ref{la=sup} can be restated as
\begin{equation}
\max _{E \in {\cal F}_n} \la(E)
= \max_{\mu} \max_{z_s}
    \sum_{s, t \in \NNN} |\lsp z_s, z_t \rsp |\mu_s \mu_t.
\label{john}
\end{equation}
where the double maximum is taken over all discrete
probability measures $\mu = (\mu_s)_{s \in \NNN}$
and all sets of unit vectors
$(z_s)_{s \in \NNN} \subset S^{n-1}(\Kn{})$ such that
$$
Id_{\KKn{n}} = n \sum_{s \in \NNN} \mu_s z_s \otimes z_s.
$$

\exm
In $\Rn{4}$, consider the following 10 vectors of the form
$$
x_{s} = \frac{1}{\sqrt{12}}
 \left( \begin{array}{c}
3 \\ -1 \\ -1\\ -1
\end{array}\right),
\qquad\quad
x_{s} = \frac{1}{\sqrt{2}}
 \left( \begin{array}{c}
\sin \al \\ \sin \al \\ -\cos \al \\ -\cos \al 2
\end{array}\right),
$$
permuting the 3 to all places in the first type of vectors
($1 \le s \le 4$) and permuting the two $-\cos \al$
in the second type of vectors
($5 \le s \le 10$).
Set $a = -  \sin 2\al  + 1/2$.
One checks that
$$
4 \Bigl( \sum_{s=1}^4 (a/ 2 (1+2a)) x_s \otimes x_s
     +  \sum_{s=5}^{10} (1/6 (1+2a)) x_s \otimes x_s  \Bigr) = Id_{\RRn{4}}.
$$
Hence, letting $ \mu_s = a/ 2 (1+2a) $ for $1 \le s \le 4$ and
$\mu_s = 1/6 (1+2a) $ for $ 5 \le s \le 10$, we see that
the $x_s$'s and $\mu_s$'s
satisfy the constraints in (\ref{john}).
The scalar products
$ |\lsp x_{s}, x_{t} \rsp|$
satisfy the following:
for $1 \le s \ne t \le 4$
they are equal to $1/3$;
for $5 \le s \ne t \le 10$ they take  two values,
$(1- \sin 2 \al)/2$
appears 24 times and $|\sin 2 \al|/2$
appears 6 times;
for $1 \le s \le 4$ and $ 5 \le t \le 10 $, they are
equal to $(1/ \sqrt 6)|\sin \al + \cos \al|$.

The maximum of the function
$$
\sum_{s, t =1}^{10} |\lsp x_{s}, x_{t} \rsp|\mu_s \mu_t
$$
is equal to $1.8494$ and it is attained for
$\al = 1.4592$.

In $\Rn{4}$, 10 equiangular vectors do note exist.
By Theorem~\ref{main} and (\ref{john}), the maximal
projection constant $\la = \sup \la (E_4)$, for
4-dimensional real spaces $E_4$ satisfies
$$
1.8494 \le \la < (2 + 3 \sqrt 6) / 5 \sim 1.8697.
$$

\medskip

The known explicit examples of equiangular lines
allow to write down the extremal
norms in the cases mentioned above, using (\ref{norm}).

\begin{center}
{\renewcommand{\arraystretch}{2.0}
\newlength{\bxwidth}
\small
\setbox2=\vtop{\hbox{dodeca-}\hbox{hedron}}
\setbox3=\hbox{%
   $\max\left(
      \max_{1\le i < j \le 7} |\al_i + \al_j|,
      \max_{1\le j \le 7} |\sum^7_{i=1,i\ne j} \al_i|
   \right)$%
}
\setlength{\bxwidth}{\wd3}
\begin{tabular}{|c|p{\bxwidth}|c|l|}
\hline
$\Kn{n}$ & \centering{$\|(\al_j)^n_1\|$} & $\la(X)$ & \\
\hline\hline
$\Rn{2}$
  & $
    \max(|2\al_1|,
    |\al_1 - \sqrt{3}\al_2|, |\al_1 + \sqrt{3}\al_2|)$
  & $ 4/3 $
  & hexagon \\
$\Rn{3}$
&
{\raggedright{%
  $
    \max_\pm(|\tau\al_1 \pm \sigma\al_2|,
    |\tau\al_2 \pm \sigma\al_3|, |\tau\al_3 \pm \sigma\al_1|)$\\
    where
     $\tau := \sqrt{{(\sqrt{5} + 1)}/{2}}$,
     $\sigma := \sqrt{{(\sqrt{5} - 1)}/{2}} $}}
  &
   \raisebox{-2.3ex}
   {$ \displaystyle \frac{\sqrt{5} + 1}{2}$}
&  \box2
   \\
$\Rn{7}$
&
\box3
&
 $ 5/2$
&
  \\
$\Rn{23}$
&
the norm is connected to points in the Leech lattice
&
$ 14/3 $
&
\\
$\Cn{2}$
&
  $
  \max(|\sqrt{3}\al_1+\al_2|,
       |\al_1+\sqrt{3}\al_2|,
       |\al_1+i\al_2|,
       |\al_1-i\al_2|)
  $
&
$\displaystyle \frac{1+\sqrt{3}}{2}$
&
\\
$\Cn{3}$
&
{\raggedright{%
  $
  \max_{j=1,2,3}
      (|\al_j-\al_{j+1}|,
       |\al_j-\omega\al_{j+1}|,
       |\al_j-\omega^2\al_{j+1}|)
  $\\
  \rule[-3mm]{0pt}{3mm}where $\omega := \exp(2\pi/3)$ and $\al_4 := \al_1$}}
&
5/3
&
\\
\hline
\end{tabular}
}
\end{center}


%
\end{document}